\documentclass[12pt]{amsart}
\usepackage{verbatim}
\usepackage{amssymb}
\usepackage{enumerate}
\newif\ifdeveloping

\newtheorem{theorem}{Theorem}[section]
\newtheorem{lemma}[theorem]{Lemma}
\newtheorem{corollary}[theorem]{Corollary}
\newtheorem{fact}[theorem]{Fact}

\newtheorem{problem}[theorem]{Problem}
\newtheorem{qtheorem}{Theorem}

\theoremstyle{definition}

\newtheorem{definition}[theorem]{Definition}

\theoremstyle{remark}
  
{}

\ifdeveloping
\usepackage[notref,notcite]{showkeys}
\fi

\newcommand{\prtime}{{\count0=\time\divide\count0 by 60
\count1=-\count0\multiply\count1 by 60
\advance\count1 by \time
\the\count0:\the\count1}
}

\def\myheads#1;#2;{
\pagestyle{myheadings}
\markboth{{\sc\hfill #1\hfill\protect\makebox[0cm][r]{\rm\today; \prtime}}}
{{\sc\protect\makebox[0cm][l]{\rm\today;\ \prtime}\hfill #2\hfill}}
\thispagestyle{myheadings}
}

\newcommand{\dcal}{{\mathcal D}}
\newcommand{\ecal}{{\mathcal E}}
\newcommand{\fcal}{{\mathcal F}}
\newcommand{\hcal}{{\mathcal H}}

\newcommand{\mcal}{{\mathcal M}}

\newcommand{\ult}[1]{\mathfrak{#1}}
\newcommand{\setm}{\setminus}
\newcommand{\empt}{\emptyset}
\newcommand{\subs}{\subset}
\newcommand{\dom}{\operatorname{dom}}

\def\<{\left\langle}
\def\>{\right\rangle}
\def\cf{\operatorname{cf}}

\def\oo{\omega_1}
\def\br#1;#2;{\bigl[ {#1} \bigr]^ {#2} }
\def\to{\longrightarrow}
\theoremstyle{plain}

\def\cl#1{\overline{#1}}
\def\cel{\operatorname{c}}
\def\hal{\operatorname{H}}

\def\Seq{\operatorname{Seq}}

\newcommand{\doe}[1]{\operatorname{\mathcal M}(#1)}

\newcommand{\conc}{^\frown}
\newcommand{\monotone}{monotonically}

\begin{document}

\author[I. Juh\'asz]{Istv\'an Juh\'asz}
\address{Alfr{\'e}d R{\'e}nyi Institute of Mathematics}
\email{juhasz@renyi.hu}

\author[L. Soukup]{
Lajos Soukup }
\address{
Alfr{\'e}d R{\'e}nyi Institute of Mathematics }
\email{soukup@renyi.hu}

\author[Z. Szentmikl\'ossy]{
Zolt\'an Szentmikl\'ossy}
\address{E\"otv\"os University of Budapest}
\email{zoli@renyi.hu}

\title{Resolvability and monotone normality}
\thanks{The preparation of this paper was supported by OTKA grant no. 61600}

\begin{abstract}
A space $X$ is said to be $\kappa$-resolvable (resp. almost
$\kappa$-resolvable) if it contains $\kappa$ dense sets that are
pairwise disjoint (resp. almost disjoint over the ideal of nowhere
dense subsets). $X$ is maximally resolvable iff it is
$\Delta(X)$-resolvable, where $\Delta(X) = \min\{ |G| : G \ne
\emptyset \mbox{ open}\}.$

We show that every crowded monotonically normal (in short: MN) space
is $\omega$-resolvable and almost $\mu$-resolvable, where $\mu =
\min\{ 2^{\omega}, \,\omega_2 \}$. On the other hand, if $\kappa$ is
a measurable cardinal then there is a MN space $X$ with $\Delta(X) =
\kappa$ such that no subspace of $X$ is $\omega_1$-resolvable.

Any MN space of cardinality $< \aleph_\omega$ is maximally
resolvable. But from a supercompact cardinal we obtain the
consistency of the existence of a MN space $X$ with $|X| = \Delta(X)
= \aleph_{\omega}$ such that no subspace of $X$ is
$\omega_2$-resolvable.
\end{abstract}

\maketitle
\ifdeveloping
\myheads{Resolvability};{Resolvability};
\fi

\section{$\omega$-resolvability}
For a topological space $X$ we denote by $\dcal(X)$ the family of
all dense subsets of $X$ and by $\mathcal{N}(X)$ the ideal of all
nowhere dense sets in $X$. Given a cardinal ${\kappa}>1$, the space
$X$ is called {\em ${\kappa}$-resolvable} iff it contains ${\kappa}$
many disjoint dense subsets. We say that $X$ is {\em almost
$\kappa$-resolvable} if there are $\kappa$ many dense sets whose
pairwise intersections are nowhere dense, that is we have
$\{D_\alpha : \alpha < \kappa\} \subs \mathcal{D}(X)$ such that
$D_\alpha \cap D_\beta \in \mathcal{N}(X)$ if $\alpha \ne \beta$.
$X$ is maximally resolvable iff it is $\Delta(X)$-resolvable, where
$\Delta(X) = \min\{ |G| : G \ne \emptyset \mbox{ open}\}$ is called
the dispersion character of $X$. Finally, if $X$ is {\em not}
$\kappa$-resolvable then it is also called $\kappa$-{\em
irresolvable}.

The following simple but useful fact, in the case of
$\kappa$-resolvability, was observed by El'kin in \cite{El}.

\begin{lemma}\label{lm:elkin}
A space $X$ is ${\kappa}$-resolvable (almost $\kappa$-resolvable)
iff every nonempty open set in $X$ includes a nonempty (and open)
${\kappa}$-resolvable (almost $\kappa$-resolvable) subset.
\end{lemma}

The aim of this paper is to present several results concerning the
(almost) resolvability properties of monotonically normal spaces.
Let us therefore recall their definition. For any topological space
$X$ we write $$\doe {X}=\bigl\{\<x,U\>\in X\times {\tau(X)}:x\in
U\bigr\}.$$ The elements of $\mathcal{M}(X)$ will be referred to as
{\em marked} open sets. The space $X$ is called {\em \monotone\
normal} iff it is $T_1$ and it admits a {\em monotone normality
operator}, that is a function $H:\doe {X}\to {\tau(X)}$ such that
\begin{enumerate}[(1)]
\item $x\in \hal(x,U)\subs U$ for each $\<x,U\>\in \doe {X}$,
\item  if
$\hal(x,U)\cap \hal(y,V) \ne \empt$ then $x \in V$ or $y \in U$.
\end{enumerate}

\def\nea{neighbourhood assignment }

We call a set $D$ in a space $X$ {\em strongly discrete} if the
points in $D$ may be separated by pairwise disjoint neighborhoods.
It is well-known that in a \monotone\ normal space any discrete
subset is strongly discrete. On the other hand, in \cite{DTTW} it
was proved that every non-isolated point of a \monotone\ normal
space is the accumulation point of a discrete subspace.
Consequently, one obtains the following result.

\begin{theorem}[{\cite{DTTW}}]\label{tm:mond}
In a \monotone\ normal space any non-isolated point is the
accumulation point of some strongly discrete set.
\end{theorem}
Let us say that a space $X$ is SD if it has the property described
in theorem \ref{tm:mond}, that is every non-isolated point of $X$ is
the accumulation point of some strongly discrete set.
\begin{theorem}\label{tm:sdis}
Any crowded SD space $X$  is ${\omega}$-resolvable.
\end{theorem}

\begin{proof}
The SD property is clearly hereditary for open subspaces. Hence, by
lemma \ref{lm:elkin}, it suffices to prove that $X$ includes an
${\omega}$-resolvable subspace.

First we show that for every strongly discrete $D\subs X$ there is a
strongly discrete $E\subs X\setm \overline{D}$ such that $D\subs
\overline E$. Indeed, fix a  \nea $U_d$ on $D$ that separates $D$
and for each $d\in D$ pick a strongly discrete set $E_d\subs X\setm
\{d\}$ with $d\in \overline{E_D}$. Then $E=\bigcup_{d\in D}(E_d\cap
U_d)$ is clearly as claimed.

Now pick an arbitrary point $x\in X$ and set $E_0=\{x\}$. Using the
above claim, for each $n<{\omega}$ we can inductively define a
strongly discrete set $E_{n+1}\subs X\setm \overline {E_n}$ such
that $E_n\subs \overline{E_{n+1}}$. Since $\cup_{i\le n}E_i\subs
\overline{E_n}$, the sets $\{E_n : n<{\omega}\}$ are pairwise
disjoint. Let us finally set $E = \bigcup\{E_n:n<{\omega}\}.$ It is
clear from our construction that if $I \subs \omega$ is infinite
then  $\bigcup\{E_n:n \in I\}$ is dense in $E$, so the subspace $E$
of $X$ is obviously $\omega$-resolvable.
\end{proof}

\begin{corollary}\label{tm:mon-norm}
Every crowded \monotone\ normal space is ${\omega}$-resol\-vable.
\end{corollary}
\bigskip

\section{$\hal$-sequences and almost resolvability}

The main result of the previous section, namely that (crowded)
\monotone\ normal spaces are $\omega$-resolvable, used very little
of the particular structure provided by monotone normality. In this
section we shall describe a procedure on \monotone\ normal spaces
that is quite specific in this respect and so, not surprisingly, it
leads to some stronger (almost) resolvability results. This
procedure had been originated (in a different form) by S. Williams
and H. Zhou in \cite{WZ}.

\begin{definition}
Let $\hal$ be a monotone normality operator on a space $X$. A family
$\mathcal{E} \subs \mathcal{M}(X)$ of marked open sets is said to be
$\hal$-{\em disjoint} if for any two members $\langle x,U \rangle,
\langle y,V \rangle$ of $\mathcal{E}$ we have $\hal(x,U) \cap
\hal(y,V) = \empt$. Clearly, if $\mathcal{E}$ is $\hal$-disjoint
then $D(\mathcal{E}) = \{ x : \exists U \mbox{ with } \langle x,U
\rangle \in \mathcal{E}\}$ is (strongly) discrete.

By Zorn's lemma, for every open set $G$ in $X$ we can fix a {\em
maximal}\\ $\hal$-disjoint family $\mathcal{E}(G) \subs
\mathcal{M}(G)$ with the additional property that $\overline{U}
\subs G$ whenever $\langle x,U \rangle \in \mathcal{E}(G)$. The
maximality of $\mathcal{E}(G)$ implies that $$\bigcup \{ \hal(x,U) :
\langle x,U \rangle \in \mathcal{E}(G)  \}$$ is a dense open subset
of $G$.
\end{definition}

With the help of the above defined operator $\mathcal{E}(G)$ we may
now describe our basic procedure as follows.

\begin{definition}\label{df:hs}
A sequence $\langle \mathcal{E}_\alpha : \alpha < \delta \rangle$ is
called a completed $\hal$-sequence of $X$ iff

\begin{enumerate}
\item $\ecal_0=\mathcal{E}(X)$,
\item for each $\alpha < \delta$ we have $$\mathcal{E}_{\alpha+1} =
\bigcup \big\{ \mathcal{E}\big(\hal(x,U) \backslash \{ x \}\big) :
\langle x,U \rangle \in \mathcal{E}_\alpha \big\},$$
\item if $\alpha < \delta$ is limit then the family $$\mathcal{W}_\alpha = \{ W \in \tau(X)
: \forall \beta < \alpha\,\, \exists \langle x,U \rangle \in
\mathcal{E}_\beta \mbox{ with } W \subs U \}$$ is a $\pi$-base in
$X$ (or, equivalently, its union $\cup \mathcal{W}_\alpha$ is dense
in $X$) and $\mathcal{E}_\alpha$ is a maximal $\hal$-disjoint
collection of marked open sets $\langle y,V \rangle$ with $V \in
\mathcal{W}_\alpha$,
\item
$$\mathcal{W}_\delta = \{ W \in \tau(X)
: \forall \beta < \delta\,\, \exists \langle x,U \rangle \in
\mathcal{E}_\beta \mbox{ with } W \subs U \}$$ is {\em not} a
$\pi$-base in $X$.

\end{enumerate}

\end{definition}

The reader may convince himself by a straight-forward transfinite
induction that the following fact is valid.

\begin{fact}
Every crowded \monotone\ normal space $X$, with monotone normality
operator $\hal$, admits a completed $\hal$-sequence $ \langle
\mathcal{E}_\alpha : \alpha < \delta \rangle$ where $\delta$ is
necessarily a limit ordinal.
\end{fact}

We now fix some notation concerning a given completed
$\hal$-sequence $ \langle \mathcal{E}_\alpha : \alpha < \delta
\rangle$ of $X$. For any ordinal $\alpha < \delta$ we put
$D_{\alpha}= D(\mathcal{E}_\alpha)$ and
$H_{\alpha}=\bigcup\{\hal(x,U):\<x,U\>\in \ecal_{\alpha}\}$. It is
clear from our definitions that each $H_\alpha$ is dense open in
$X$, moreover $\beta < \alpha < \delta$ implies that $H_\beta
\supset H_\alpha$ and $D_\beta \cap H_\alpha = \empt$. If $I \subs
\delta$ is a set of ordinals we write
$D[I]=\bigcup\{D_{\alpha}:{\alpha}\in I\}$. Finally, we set $V = X
\setminus \overline{\cup \mathcal{W}_\delta}$, then $V$ is a
non-empty open set in $X$.

\begin{lemma}\label{lm:nwd}
If $I$ is bounded in ${\delta}$ then $D[I]$ is nowhere dense in $X$.
However, if $I$ is unbounded in ${\delta}$ then $D[I]$ is dense in
$V$, that is \\ $V \subset \overline{D[I]}$.
\end{lemma}

\begin{proof}
The first part is obvious because $I \subs \alpha < \delta$ implies
$D[I] \cap H_\alpha = \empt$.

Assume now that $I$ is cofinal in $\delta$ but, arguing indirectly,
for some $G\in {\tau}^*(V)$ we have $G\cap D[I]=\empt$. Pick any
point $z\in G$,
we claim that then, for all ${\alpha}<{\delta}$ and $\<x,U\>\in
\ecal_{\alpha}$, $\,\,\hal(x,U)\cap \hal(z,G) \ne \empt$ implies
$z\in \hal(x,U)$.

Indeed, if ${\beta}\in ({\alpha},{\delta})\cap I$ then there is
$\<x',U'\>\in \ecal_{\beta}$ with $$\hal(x',U')\cap \hal(x,U)\cap
\hal(z,G) \ne \empt$$ because $H_\beta$ is dense in $X$. It follows
that $U'\subs \hal(x,U)$, hence $x'\notin G$ as $x' \in D_\beta$ and
$G \cap D_\beta = \empt$. But then $\hal(x',U')\cap\hal(z,G)\ne
\empt$ implies $z\in U'\subs \hal(x,U)$.

The sets $\{\hal(x,U):\<x,U\>\in\ecal_{\alpha}\}$ being pairwise
disjoint, it follows that for each ${\alpha}<{\delta}$ there is
exactly one $\<x_{\alpha},U_{\alpha}\>\in \ecal_\alpha$ such that \\
$\hal(x_\alpha, U_{\alpha}) \cap \hal(z,G) \ne\empt$. But then
$\hal(z,G) \subs \cl{\hal(x_\alpha, U_\alpha)}
 \subs \cl {U_{\alpha}}$ whenever ${\alpha}<{\delta}$, consequently $$\hal(z,G) \subs
\overline{U_{{\alpha}+1}} \subs U_{\alpha}$$ as well. This, however,
would imply $\hal(z, G) \in \mathcal{W}_{\delta}$, contradicting
that $\hal(z, G) \subs G \subs V.$
\end{proof}

We may now give the main result of this section.

\begin{theorem}\label{tm:ma}
Any crowded \monotone\ normal space $X$ is almost
$\min(\mathfrak{c}, \omega_2)$-resolvable. So $X$ is almost
$\omega_1$-resolvable, and even almost $\omega_2$-resolvable if the
continuum hypothesis (CH) fails.
\end{theorem}

\begin{proof}
By lemma \ref{lm:elkin} it suffices to show that some non-empty open
$V \subs X$ satisfies this property. To see this, let us consider a
completed $\hal$-sequence $\langle \mathcal{E}_\alpha : \alpha <
\delta \rangle$ of $X$. Let $I$ be a cofinal subset of ${\delta}$ of
order type $\cf({\delta})$ and $\{I_{\zeta}:{\zeta}<{\mu}\}$ be an
almost disjoint subfamily of $[I]^{\cf(\delta)}$, where
${\mu}=2^{\omega} = \mathfrak{c}$ if $\cf({\delta})={\omega}$ and
${\mu}=\cf({\delta})^+ \ge \omega_2$ if $\cf({\delta})>{\omega}$.
Then the family $\{D[I_{\zeta}]:{\zeta}<{\mu}\}$ witnesses that $V$
is almost ${\mu}$-resolvable.
\end{proof}

Since almost $\omega$-resolvability is clearly equivalent
$\omega$-resolvability, theorem \ref{tm:ma} provides us a new proof
of \ref{tm:mon-norm}.

\section{spaces from trees and ultrafilters}

The prime examples of \monotone\ normal spaces are metric and
ordered spaces that are all known to be maximally resolvable.
Compared to this the results of the two preceding sections seem
rather modest. The main aim of this section is to show that, at
least modulo some large cardinals, nothing stronger than
$\omega$-resolvability can be expected of a \monotone\ normal space
$X$, even if the dispersion character $\Delta(X)$ is large. The
examples that show this have actually been around but, as far as we
know, the fact that they are \monotone\ normal has not been noticed.

The underlying set of such a space is an {\em everywhere infinitely
branching} tree $\langle T,< \rangle$. This simply means that for
each $t \in T$ the set $S_t$ of all immediate successors of $t$ in
$T$ is infinite. The height of such a tree is obviously a limit
ordinal. (In fact, nothing is lost if we only consider trees of
height $\omega$.) By a {\em filtration} on $T$ we mean a map $F$
with domain $T$ that assigns to every $t \in T$ a filter $F(t)$ on
$S_t$ such that every cofinite subset of $S_t$ belongs to $F(t)$
(that is, $F(t)$ extends the Fr\'echet filter on $S_t$).

\begin{definition}
Assume that $F$ is a filtration on an everywhere infinitely
branching tree $\langle T,< \rangle$. A topology $\tau_F$  is then
defined on $T$ by $$\tau_F = \{ V \subs T : \forall t \in V\,\,
\big(V \cap S_t \in F(t) \big) \},$$ and the space $\langle T,\tau_F
\rangle$ is denoted by $X(F)$.
\end{definition}

\begin{theorem}\label{tm:tr}
Let $F$ be a filtration on an everywhere infinitely branching tree
$\langle T,< \rangle$. Then the space $X(F)$ is \monotone\ normal.
\end{theorem}

\begin{proof}
That $\tau_F$ is indeed a topology that satisfies the $T_1$
separation axiom is obvious and well-known. The novelty is in
showing that $X(F)$ is \monotone\ normal.

To this end we define $\hal(s,V)$ for $s \in V \in \tau_F$ as
follows: $$\hal(s,V) = \{ t \in V : s \le t \mbox{ and } [s,t] \subs
V \}.$$ Of course, here $[s,t] = \{ r : s \le r \le t \}$. Clearly,
$\hal(s,V) \in \tau_F$ and $s \in \hal(s,V) \subs V.$

Next, assume that $t \in \hal(s_1,V_1) \cap \hal(s_2,V_2).$ Then
$s_1, s_2 \le t$ implies that $s_1$ and $s_2$ are comparable, say
$s_1 \le s_2$. But then we have \\ $s_2 \in [s_1,t] \subs V_1$,
consequently $\hal$ is indeed a monotone normality operator on
$X(F).$
\end{proof}

Of special interest are those filtrations $F$ for which $F(t)$ is a
(free) ultrafilter on $S_t$ for all $t \in T$. Such an $F$ will be
called an {\em ultrafiltration}. In this case we have a convenient
way to determine the closures of sets in the space $X(F)$ that will
be put to good use later.

\begin{definition}
 For every set $A
\subs T$ we define $$C(A) = A \cup \{ t \in T : S_t \cap A \in
F(t)\}.$$ Then by transfinite recursion we define $C^\alpha(A)$ for
all ordinals $\alpha$ by $C^{\alpha+1}(A) = C(C^\alpha(A))$ for
successors and $C^\alpha(A) = \cup \{ C^\beta : \beta < \alpha\}$
for $\alpha$ limit.
\end{definition}

\begin{lemma}\label{lm:cl}
Let $F$ be an ultrafiltration on the tree $T$. Then a set $B \subs
T$ is closed in $X(F)$ iff $B = C(B)$. Consequently, for any subset
$A \subs T$ there is an ordinal $\alpha < |T|^+$ with $\overline{A}
= C^\alpha(A).$
\end{lemma}

\begin{proof}
First, if $B = C(B)$ then for each $t \in T \backslash B$ we have
$S_t \cap B \notin F(t)$, hence $S_t \backslash B \in F(t)$ because
$F(t)$ is an ultrafilter. Then $T \backslash B$ is open by the
definition of $\tau_F$, hence $B$ is closed. Conversely, if $B$ is
closed in $X(F)$ then for each $t \in T \backslash B$ we have $S_t
\backslash B \in F(t)$, hence $S_t \cap B \notin F(t)$, that is $t
\notin C(B).$ But this means that $B = C(B)$.

Next, $C(A) \subs \overline{A}$ is obvious, and then by induction we
get $C^\alpha(A) \subs \overline{A}$ for all $\alpha$. But for some
$\alpha < |T|^+$ we must have $C(C^\alpha(A)) = C^\alpha(A)$, and
then $\overline{A} = C^\alpha(A)$ for $C^\alpha(A)$ is closed by the
above.
\end{proof}

Let  $\ult u$ be an ultrafilter on a set $I$ and $\lambda$ be a cardinal.
$\ult u$ is said to be {\em $\lambda$-descendingly complete } iff
$\bigcap\{X_\xi:\xi<\lambda\}\in {\ult u}$ for each
decreasing sequence
 $\{X_\xi:\xi<\lambda\}\subs \ult u$.
The ultrafilter ${\ult u}$ is called {\em $\lambda$-descendingly
incomplete} iff it is not  $\lambda$-descendingly complete. For
example,  $\ult u$ is countably complete exactly if it is
$\omega$-descendingly complete.

We shall need the following old result of Kunen and Prikry in our
next irresolvability theorem for spaces obtained from certain
ultrafiltrations.

\begin{qtheorem}[Kunen, Prikry, \cite{KP}]\label{tm:KP}
If $\lambda$ is a regular cardinal and $\ult u$ is a
$\lambda$-descendingly complete ultrafilter (on any set) then $\ult
u$ is also $\lambda^+$-descendingly complete.
\end{qtheorem}

\begin{theorem}\label{tm:subad-resolv}
Assume that $F$ is an ultrafiltration on $T$ and ${\lambda}$ is a
regular cardinal such that $F(t)$ is ${\lambda}$-descendingly
complete for all $t \in T$. Then the space $X(F)$ is hereditarily
${\lambda}^+$-irresolvable.
\end{theorem}

\begin{proof}
First we show that for every set $A \subs T$ we have
{$\overline{A}=C^{\lambda}(A)$}. By lemma \ref{lm:cl} it suffices to
show that $C(C^{\lambda}(A))=C^{\lambda}(A)$.

Assume, indirectly, that $t\in C(C^{\lambda}(A)) \backslash
C^{\lambda}(A)$, then we must have $C^\lambda(A) \cap S_t \in F(t)$.
But
$$C^\lambda(A) \cap S_t = \bigcup_{\alpha < \lambda}C^\alpha(A) \cap S_t$$
where the right-hand side is an increasing union, hence there is an
$\alpha < \lambda$ with $C^\alpha(A) \cap S_t \in F(t)$ because
$F(t)$ is $\lambda$-descendingly complete. This implies that $t \in
C^{\alpha+1}(A) \subs C^\lambda(A)$, a contradiction.

\newcommand{\ord}{\operatorname{ord}}
Let us now consider an {\em indexed} family of sets $\fcal = \{ F_i
: i \in I \}$. We are going to use the following notation:
\begin{displaymath}
\ord (x,\fcal)=|\{i\in I :x\in F_i\}|
\end{displaymath}
and
\begin{displaymath}
\ord (\fcal)=\sup\{\ord(x,\fcal):x\in\cup_{i \in I}F_i\}.
\end{displaymath}

Instead of the statement of the theorem we shall prove the following
much stronger claim.
\begin{lemma}\label{lm:ord}
If  $\dcal = \{D_i : i \in I\}$ is any indexed family of subsets of
$T$ with $\ord (\dcal) \le {\lambda}$ then $\ord (\{\overline{D_i}:
i \in I\}) \le {\lambda}$ as well.
\end{lemma}

\begin{proof}
We shall prove, by induction on ${\alpha}\le {\lambda}$, that
$\ord(\mathcal{D}^\alpha) \le \lambda$ where
\begin{displaymath}
\mathcal{D}^\alpha = \{C^{\alpha}(D_i):i\in I\} .
\end{displaymath}
\noindent We first show that $\ord(\mathcal{D}^1) \le \lambda$, this
will clearly take care of all the successor steps.

Assume, indirectly, that $\ord(t,\mathcal{D}^1) \ge \lambda^+$ for
some $t \in T$, then we may find a set $J \in [I]^{\lambda^+}$  such
that $t \in C(D_j) \backslash D_j$, hence $D_j \cap S_t \in F(t)$,
for each $j \in J$.

By the theorem of Kunen and Prikry the ultrafilter $F(t)$ is also
${\lambda}^+$-descendingly complete. Consequently, using a standard
argument, one can show that there is an $L\in \br J;{\lambda}^+;$
such that
$$\bigcap \{D_j \cap S_t : j \in L\} \ne \emptyset.$$ But this
clearly contradicts $\ord(\mathcal{D}) \le \lambda$.

Next assume that ${\alpha} \le \lambda$ is a limit ordinal and the
inductive hypothesis holds for all $\beta < \alpha$. But now for
each index $i \in I$ we have $C^\alpha(D_i) = \bigcup_{\beta <
\alpha}C^\beta(D_i)$, hence
$$\ord(t,\mathcal{D}^\alpha) \le \sum_{\beta < \alpha} \ord
(t,\mathcal{D}^\beta) \le |\alpha|\cdot\lambda =\lambda$$ whenever
$t \in T$, and so $\ord(\mathcal{D}^\alpha) \le \lambda$.
\end{proof}

It follows immediately from lemma \ref{lm:ord} that if $\{A_i : i
\in \lambda^+\}$ are pairwise disjoint non-empty subsets of $T$ then
the closures $\overline{A_i}$ cannot all be the same and so no
subspace of $X(F)$ can be $\lambda^+$-resolvable.
\end{proof}

\begin{corollary}\label{co:me}
If $F$ is an ultrafiltration on $T$ such that $F(t)$ is countably
complete for each $t \in T$ then $X(F)$ is $\omega$-resolvable but
hereditarily $\omega_1$-irresolvable. In particular, if $\kappa$ is
a measurable cardinal then there is a \monotone\ normal space $X$
with $|X| = \Delta(X) = \kappa$ that is hereditarily
$\omega_1$-irresolvable.
\end{corollary}

The question if $\omega$-resolvable spaces are also maximally
resolvable was raised a long time ago by Ceder and Pearson in
\cite{CP}, and has just recently been settled completely
 in \cite{JSSz} (negatively). Corollary \ref{co:me} yields a
\monotone\ normal counterexample to this problem, from a measurable
cardinal. Another counterexample from a measurable cardinal was
found by Eckertson in \cite{E}, however, that example is not
\monotone\ normal. We present two arguments to show this. First,
Eckertson's example contains a crowded irresolvable subspace, hence
it cannot be \monotone\ normal by corollary \ref{tm:mon-norm}.

The second argument is based on our following observation that may
have some independent interest. First we need some  notation. If
${\kappa}\le {\lambda}$ are cardinals we let
${\tau}^{\lambda}_{\kappa}$ denote the $\,\,< \kappa$ box product
topology on $2^\lambda$ (generated by the base $\{[f]:f\in
Fn({\lambda},2;{\kappa}\}$, where $[f]=\{x\in 2^\lambda : f\subs
x\}$), moreover we set $\mathbb{C}_{\lambda,\kappa} =
\<2^{\lambda},{\tau}^{\lambda}_{\kappa}\>.$

\begin{theorem}\label{tm:box}
If ${\kappa}^{<{\kappa}}={\kappa}<{\lambda}$ then no dense subspace
of $\mathbb{C}_{\lambda,\kappa}$  is \monotone\ normal.
\end{theorem}

\begin{proof}[Proof of \ref{tm:box}]

\def\doe{{\dot\tau}}
\newcommand{\supp}{\operatorname{supp}}

Let $X$ be dense in $\mathbb{C}_{\lambda,\kappa}$ and $\theta$ be a
large enough regular cardinal. Let $\mcal$ be an elementary submodel
of  $\<\hcal(\theta),\in, \prec\>$ (where $\hcal(\theta)$ is the
family of sets hereditarily of size $<\theta$ and $\prec$ is a
well-ordering of $\hcal(\theta)$) such that $|\mcal|={\kappa}$ and
$\br \mcal;<{\kappa};\subs \mcal$, moreover $X,
{\kappa},{\lambda}\in \mcal$. Note that then $Fn(\br \mcal\cap
{\lambda};<{\kappa};,2;{\kappa})\subs \mcal$ as well.

Assume that $X$ is \monotone\ normal and let $\hal\in \mcal $ be  a
monotone normality operator on $X$. We can assume that $\hal(x,[s]
\cap X)$ is the trace on $X$ of a basic open set for each basic open
set $[s]$.

Let $I=\mcal\cap {\lambda}$ and pick ${\alpha}\in {\lambda}\setm I$.
$\fcal=\{f\restriction I :f\in \mcal\cap X\}$ is clearly dense in
the subspace $2^I$ of $\mathbb{C}_{\lambda,\kappa}$. Let
$\fcal_i=\{f\restriction I:f\in X\cap \mcal\land f({\alpha})=i\}$
for $i\in 2$ then $\fcal=\fcal_0\cup\fcal_1$ so there is $i\in 2$
and $s\in Fn(I,2;{\kappa})$ such that $\fcal_i$ is dense in $2^I\cap
[s] \cap X$.

Let $b=s\cup\{\<{\alpha},1-i\>\}$ and pick $x\in X\cap [b]$. Next,
let $\hal(x,[b] \cap X) = [b'] \cap X$  and $b''=b'\restriction I$.
Fix $b'''\in Fn(I,2;{\kappa})$ such that $b'''\supset b''$ and
$x\notin [b''']$. Since $\fcal_i$ is dense in $2^I\cap [s] \cap X$
we can pick $y\in X\cap\mcal\cap [b''']$ such that $y({\alpha})=i$.
Let $[u] \cap X =\hal(y,[b'''] \cap X)$. Then $\dom u\subs I$
because $\hal,b''',y\in\mcal$.

Since $x\notin [b''']$ and $y\notin [b]$ it follows that
$\hal(x,[b])\cap \hal(y,[b'''])=[u]\cap [b']\cap X=\empt$. However
$\supp u\subs I$ and $u\supset b'''\supset b''=b'\restriction I$, so
$u$ and $b'$ are compatible  functions of size $<{\kappa}$, i.e.
$[u]\cap [b']$ is a nonempty open set in $\<2^{\lambda},
{\tau}^{\lambda}_{\kappa}\>$. Since $X$  is dense we have $[u]\cap
[b']\cap X\ne \empt$, a contradiction.
\end{proof}

Now, Eckertson's example obtained from a measurable cardinal
$\kappa$ contains a subspace homeomorphic to a dense subspace of
$\mathbb{C}_{2^\kappa,\kappa}\,$, hence it cannot be \monotone\
normal by theorem \ref{tm:box} because
${\kappa}^{<{\kappa}}={\kappa}$.

Of course, we have a space like in corollary \ref{co:me} iff there
is a measurable cardinal. Also, the cardinality (and dispersion
character) of such a space is at least as large as the first
measurable. But can one have a \monotone\ normal example that is
much smaller? The answer to this question is, consistently, yes
assuming the existence of a large cardinal that is even stronger
than a measurable.

\begin{qtheorem}[Magidor, \cite{Mag2}] It is consistent from a
supercompact cardinal  that there is an $\oo$-descendingly complete
uniform ultrafilter on $\aleph_{\omega}$.
\end{qtheorem}

In \cite{Mag} a similar statement was proved for 
$\aleph_{\omega+1}$ instead of $\aleph_{\omega}$, but according to Magidor 
a slight modification of that   proof works even for  $\aleph_{\omega}$.

From this result of Magidor and from theorem \ref{tm:subad-resolv}
we immediately obtain our promised result.

\begin{corollary}\label{co:mag}
From a supercompact cardinal it is consistent to have a \monotone\
normal space $X$ with $|X|=\operatorname{\Delta}(X)=\aleph_{\omega}$
that is hereditarily ${\omega}_2$-irresolvable (hence not maximally
resolvable).
\end{corollary}
Actually, in \cite{Mag} a slightly weaker result is given in which
$\aleph_\omega$ is replaced with $\aleph_{\omega+1}$. However, in a
private communication, Magidor pointed out to us that the method of
\cite{Mag} yields the above stronger version as well.

But can we do even better and go below $\aleph_\omega$? The answer
to this question is, maybe surprisingly, negative.  We are going to
show that any \monotone\ normal space of cardinality less than
$\aleph_{\omega}$ is maximally resolvable. The proof of this result
will be based on showing that all spaces of the form $X(F)$ with $F$
an ultrafiltration on the tree $\Seq {\kappa} = \kappa^{<\omega}$ of
all finite sequences of ordinals less than $\kappa$ are maximally
resolvable provided that ${\kappa}<\aleph_{\omega}$. The first
result to this effect, for constant ultrafiltrations on $\Seq
\omega_n$, was obtained by L\'aszl\'o Heged\"us in his Master's
Thesis \cite{Heg}. Of course, by a constant ultrafiltration we mean
one for which $F(t)$ is the "same" ultrafilter for all $t \in T$.

Now, let $\kappa$ be an arbitrary infinite cardinal. A non-empty
subset $T$ of $\Seq {\kappa}$ is called a {\em subtree} of $\Seq
{\kappa}$ iff $t\restriction n\in T$ whenever $t\in T$ and $n<|t|$.
For any subset $A$ of $\Seq {\kappa}$ we shall write $\,\,\min
A\,\,$ to denote the set of all minimal elements of $A$ (with
respect to the tree ordering on $\Seq {\kappa}$, of course).

If $F$ is a filtration on $\Seq {\kappa}$ and $v\in \Seq {\kappa}$
we shall denote by $F_v$ the derived filtration on $\Seq {\kappa}$
defined by the formula $F_v(s)=F(v^\frown s)$.

Assume now that $S$ and $\{T_v:v\in \Seq {\kappa}\}$ are subtrees of
$\Seq {\kappa}$. We then define their ``sum" by
\begin{displaymath}
S\oplus \{T_v:v\in \Seq {\kappa}\}= S\cup \{v^\frown t:v\in \min
(\Seq {\kappa}\setm S) \land t\in T_v\}.
\end{displaymath}
Obviously, this sum is again a subtree of $\,\,\Seq \kappa$.

If moreover $f$ and $g = \{g_v:v\in \Seq {\kappa}\}$ are functions
with $\dom f=S$ and $\dom g_v=T_v$ then we define $f\oplus
\{g_v:v\in \Seq {\kappa}\} = f \oplus g$ by putting
$$\dom (f\oplus g)=S\oplus \{T_v:v\in \Seq {\kappa}\}$$ and
\begin{displaymath}
(f\oplus g)(x)=\left\{
\begin{array}{ll}
f(x)&\text{for $x\in S$}\\g_v(t)&\text{for $x=v\conc t$ with $v\in
\min(\Seq {\kappa}\setm S)$, $ t\in T$.}
\end{array}
\right.
\end{displaymath}

\newcommand{\wf}{well-founded}
A subtree of $\Seq {\kappa}$ is called {\em \wf} iff it does not
possess any infinite branches. Note that if $S$ and $\{T_v:v\in \Seq
{\kappa}\}$ are all {\wf} then so is $S\oplus \{T_v:v\in \Seq
{\kappa}\}$.

Now let $0 < \lambda \le{\kappa}$  be cardinals and $F$ be a
filtration on $\Seq {\kappa}$. We say that a function $f$ is {\em
${\lambda}$-good for $F$} iff $\,\,\dom f$ is a well-founded subtree
of $\Seq {\kappa}$, moreover $f[V]={\lambda}$ whenever $V$ is open
in $X(F)$ with $\empt\in V$. As an easy (but useful) illustration of
this concept we present the following result.

\begin{lemma}\label{lm:n}
For each $0 < n< {\omega}$ and for any filtration $F$ on $\kappa$
there is a function $f$ which is $n$-good  for $F$.
 \end{lemma}

\begin{proof}
Let $\dom f=\{s\in \Seq {\kappa}: |s|<n\}$ and $f(s)=|s|$.
\end{proof}

The next result shows the relevance of these concepts to
resolvability.

\begin{theorem}\label{tm:good2col}
Let $F$ be an filtration on $\Seq {\kappa}$. If there are
${\lambda}$-good functions $f_s$ for $F_s$  for all $s\in \Seq
{\kappa}$ then $X(F)$ is ${\lambda}$-resolvable.
\end{theorem}

\begin{proof}
Define the sequence of functions $\,\,g_0, g_1,\dots$ by recursion
as follows: $g_0=f_{\emptyset}$ and $g_{n+1}=g_n\oplus \{f_s:s\in
\Seq {\kappa}\}$ for $n<{\omega}$. It is easy to check that then
$g_{\omega}=\bigcup_{n<{\omega}}g_n$ maps $\Seq {\kappa}$ to
${\lambda}$, i. e. $\dom g_\omega = \Seq {\kappa}$. Indeed, if $s
\in \Seq \kappa$ with $|s| = n$ then there is a $k \le n$ with $s
\in \dom g_k$.

We show next that ${g_{{\omega}}}[V]={\lambda}$ holds for any
non-empty open set $V$ in $X(F)$. Let $n$ be such that $V \cap \dom
g_n \ne \empt$ and pick $v \in V \cap \dom g_n$. Clearly, there is
an extension $s$ of $v$ with $s\in V \cap \min (\Seq {\kappa}\setm
\dom g_n).$ Now let
$$W=\{t\in \Seq {\kappa}:s\conc t\in V\}$$ then $\emptyset \in W$ and $W$ is
open in $X(F_s)$, hence $f_s[W]={\lambda}$ because $f_s$ is
$\lambda$-good for $F_s$. But we clearly have $g_{\omega}(s\conc
t)=f_s(t)$ for all $t\in \dom f_s$, hence we have
${g_{\omega}}[V]={\lambda}$ as well.

But then $\{ g_\omega^{-1}(\alpha) : \alpha < \lambda \}$ is a
pairwise disjoint family of dense sets in $X(F)$.
\end{proof}

The following stepping-up type result will turn out to be very
useful.

\begin{lemma}\label{lm:good}
Assume that $F$ is a filtration on $\,\Seq{\kappa}\,\,$such that
$F(\emptyset)$ is ${\lambda}$-descendingly incomplete, moreover for
every cardinal ${\mu}<{\lambda}$ and every ordinal
${\alpha}<{\kappa}$ there is a $\mu$-good function
$f^{\alpha}_{\mu}$ for $F_{\<{\alpha}\>}$. Then there is a
${\lambda}$-good function $f$ for $F$.
\end{lemma}

\begin{proof}
Fix a continuously decreasing sequence
 $\{X_\xi:\xi<\lambda\}\subs F(\emptyset)$ with empty intersection. For any ordinal
$\nu < \lambda$ let us put $I_\nu = X_\nu \backslash X_{\nu+1}$,
then we clearly have ${\kappa}=\bigcup\{I_{\nu}:{\nu}<{\lambda}\}$.
For each  $0 <{\nu}<{\lambda}$ fix a map
$h_{\nu}:|{\nu}|\stackrel{onto}{\to} {\nu}$.

We now define the desired map $f$ with the following stipulations:
\begin{displaymath}
\dom f=\{\empt\}\cup\bigcup_{\nu<\lambda}\{\<{\alpha}\>^{\frown} t:
{\alpha}\in I_{\nu}\, \mbox{ and } t\in \dom
f_{|{\nu}|}^\alpha\}\,\,,
\end{displaymath}
and for $s \in \dom f$

\begin{displaymath}
f(s)=\left\{
\begin{array}{ll}
0&\text{if $s=\empt\,$,}\\h_{\nu}(f^{\alpha}_{|{\nu}|}(t))&\text{if
$s=\<{\alpha}\>\conc t$  with ${\alpha}\in I_{\nu},\,\, t\in \dom
f_{|\nu|}^\alpha\,\,$.}
\end{array}
\right.
\end{displaymath}

Clearly, $f$ is well-defined and $\dom f$ is \wf. If $V$ is open in
$X(F)$ with $\empt\in V$ then we have $V \cap S_\emptyset \in
F(\emptyset)$ and hence
\begin{displaymath}
\sup\{{\nu}: \exists {\alpha}\in I_{\nu}\ \text{with}
\<{\alpha}\>\in V \}={\lambda}.
\end{displaymath}
But $\<{\alpha}\>\in V$ and ${\alpha}\in I_{\nu}$ imply
${f^{\alpha}_{|{\nu}|}}[\{s:\<{\alpha}\>\conc s\in V\}]=|{\nu}|$ and
so $f[V] \supset {\nu}$, hence we have $f[V] = \lambda.$
\end{proof}

\begin{theorem}\label{tm:uj}
Let $F$ be a filtration on $\Seq \kappa$ and $\lambda$ be an
infinite cardinal such that $F(t)$ is $\mu$-descendingly incomplete
whenever $t \in \Seq \kappa$ and $\omega \le \mu \le \lambda$. Then
there are $\lambda$-good functions for all the derived filtrations
$F_s$ and hence $X(F)$ is $\lambda$-resolvable.
\end{theorem}

\begin{proof}
The proof goes by a straight-forward induction on $\lambda$, using
lemma \ref{lm:good} and the fact that our assumption on $F$ is
automatically valid also for all the derived filtrations $F_s$. The
starting case $\lambda = \omega$ also uses lemma \ref{lm:n}. The
last statement is immediate from theorem \ref{tm:good2col}.
\end{proof}

A uniform ultrafilter on $\kappa$ is trivially $\kappa$-descendingly
incomplete. So if $\kappa = \omega_n < \aleph_\omega$, then it
follows by $n$ repeated applications of the above mentioned result
of Kunen and Prikry that any uniform ultrafilter on $\kappa$ is
$\mu$-descendingly incomplete for all $\mu$ with $\omega \le \mu \le
\kappa.$ Thus we get from theorem \ref{tm:uj} the following result.

\begin{corollary}\label{co:uni}
Assume that $\kappa < \aleph_\omega$ and $F$ is any uniform
ultrafiltration on $\Seq \kappa$ ( i. e. F(t) is uniform for all $t
\in \Seq \kappa$). Then the space $X(F)$ is $\kappa$-resolvable.
\end{corollary}

We now recall a definition from \cite{jssz2}, see also \cite{Pa}.

\begin{definition}
Let $X$ be a space and $\mu$ be an infinite cardinal number. We say
that $x \in X$ is a $T_\mu$ point of $X$ if for every set $A \in
[X]^{<\mu}$ there is some $B \in [X \backslash A]^{<\mu}$ such that
$x \in \overline{B}$. We shall use $T_\mu(X)$ to denote the set of
all $T_\mu$ points of $X$.
\end{definition}

The following result is an easy consequence of lemma 1.3 from
\cite{jssz2}. In the particular case when $\mu$ is a successor
cardinal it follows from proposition 2.1 of \cite{Pa}.

\begin{lemma}\label{lm:pav}
If $|X| = \mu$ is a regular cardinal  and $T_{\mu}(X)$ is dense in
$X$ then $X$ is $\mu$-resolvable.
\end{lemma}

This result will enable us to transfer certain results from spaces
of the form $X(F)$, where $F$ is a uniform ultrafiltration on $\Seq
\kappa$ for some regular cardinal $\kappa$, to monotonically normal
and even more general spaces.

Let us recall from section 1 that every \monotone\ normal space is
SD. In fact, as monotone normality is a hereditary property, it is
even hereditarily SD (in short: HSD). We shall need below a property
that is strictly between SD and HSD, namely that all {\em dense}
subspaces are SD, we shall denote this property by DSD. It can be
shown that for instance the \v Cech-Stone remainder $\omega^*$ is
DSD but not HSD.

\begin{theorem}\label{tm:ekv}
Assume that $\kappa = \cf(\kappa) \ge \lambda.$ Then the following
are equivalent.
\begin{enumerate}
\item
If $X$ is a DSD space with $|X| = \Delta(X) = \kappa$ then $X$ is
$\lambda$-resolvable.

\item
If $X$ is a MN space with $|X| = \Delta(X) = \kappa$ then $X$ is
$\lambda$-resolvable.

\item
If $F$ is any uniform ultrafiltration on $\Seq \kappa$ then the
space $X(F)$ is $\lambda$-resolvable.
\end{enumerate}
\end{theorem}

\begin{proof}
Of course, only (3) $\Rightarrow$ (1) requires proof. So assume (3)
and consider a DSD space $X$ with $|X| = \Delta(X) = \kappa.$ If
$T_\kappa(X)$ is dense in $X$ then, by lemma \ref{lm:pav} $X$ is
even $\kappa$-resolvable and we are done.

Otherwise, in view of lemma \ref{lm:elkin}, we may assume that
actually $T_\kappa(X) = \emptyset.$ In this case for every point $x
\in X$ there is a set $A_x \in [X]^{<\kappa}$ such that $x \in A_x$
and for $D_x = X \backslash A_x$ no $B \in [D_x]^{<\kappa}$ has $x$
in its closure. Note that by $\Delta(X) = \kappa$ each $D_x$ is
dense in $X$.

But $X$ is DSD, hence for every $x$ there is a strongly discrete set
$S_x \subs D_x$ with $x \in \overline{S_x}$. (Note that $S \subs
D_x$ is strongly discrete in $D_x$ iff it is so in $X$ for $D_x$ is
dense.)

Next, by recursion on $|t|\,\, $, we define points $x_t$ and open
sets $U_t$ in  $X$ as follows. First pick any  point $x_\emptyset
\in X = U_\emptyset$. If $x_t \in U_t$ has been defined then fix a
one-to-one enumeration of $S_{x_t} \cap U_t = \{ x_{t \conc \alpha}
: \alpha < \kappa  \}$ and choose $\{ U_{t \conc \alpha} : \alpha <
\kappa \}$ to be pairwise disjoint open neighbourhoods of them, all
contained in $U_t$. Clearly, then the map $h : \Seq \kappa \to X$
that maps $t$ to $ h(t) =x_t$ is injective.

Next, for any $t \in \Seq \kappa$ extend the trace of the
neighbourhood filter of $x_t$ on $S_{x_t} \cap U_t$ to an
ultrafilter $u_t$ and define $F(t) = h^{-1}[u_t]$, which is an
ultrafilter on $S_t = \{ t \conc \alpha : \alpha < \kappa \}$. It
follows from our assumptions that every $F(t)$ is uniform and
therefore $X(F)$ is $\lambda$-resolvable. But the subspace topology
on $h[\Seq \kappa]$ in $X$ is clearly coarser than the $h$-image of
$\tau_F$, hence it is also $\lambda$-resolvable. By lemma
\ref{lm:elkin}, this completes our proof.
\end{proof}

\begin{corollary}
Let $X$ be any DSD space of cardinality $< \aleph_\omega$. Then $X$
is maximally resolvable. In particular, all MN spaces of size $<
\aleph_\omega$ are maximally resolvable.
\end{corollary}

\begin{proof}
Clearly, every open set $U$ in $X$ includes another open set $V$
such that $|V| = \Delta(V)$. But every open subspace of a DSD space
is again DSD, so theorem \ref{tm:ekv} and corollary \ref{co:uni}
imply that $V$ is $|V|$-resolvable. But $\Delta(X) \le |V|$, hence
each such $V$ is $\Delta(X)$-resolvable and so, in view of lemma
\ref{lm:elkin}, $X$ is maximally resolvable.
\end{proof}

We conclude by listing a few open problems that we find interesting.

\begin{problem}

\begin{enumerate}

\item

Is there a ZFC that example of a \monotone\ normal space that is not
maximally resolvable?
\item
Is it consistent to have a \monotone\ normal space $X$ of
cardinality less than the first measurable such that $\Delta(X) >
\omega$ but $X$ is not $\omega_1$-resolvable?
\item
Is every crowded monotonically normal space almost
$\mathfrak{c}$-resolvable?
\end{enumerate}
\end{problem}

\end{document}